\documentclass[12pt]{article}
\usepackage{eurosym}
\usepackage{amsmath,amssymb,amsthm,amscd,graphicx}
\usepackage{latexsym}
\usepackage{enumerate}
\usepackage{amsbsy, amsfonts, textcomp}
\usepackage{color}
\usepackage[english]{babel}
\usepackage{tikz-cd}

\newtheorem{theorem}{Theorem}[]

\theoremstyle{definition}
\newtheorem{definition}[theorem]{Definition}
\newtheorem{example}[theorem]{Example}
\theoremstyle{remark}
\newtheorem{remark}[theorem]{Remark}

\newcommand{\R}{\mathbb{R}}

\newcommand{\F}{\mathcal{F}}

\def \Ker{\operatorname{Ker}}

\def \Sing{\operatorname{Sing}}
\def \grad{\operatorname{grad}}
\numberwithin{theorem}{subsection}

\begin{document}

\begin{center}

{\Large
{\bf Tame topology and non-integrability of dynamical systems \let\thefootnote\relax\footnote{
R.
Mohseni acknowledges support of the Polish Ministry of Science and
Higher Education. Z. Hajto acknowledges support of
grant PID2019-107297GB-I00 (MICINN).}}

{\normalsize \vspace{1cm} }

{\normalsize \vspace{0.5cm} ZBIGNIEW HAJTO \\[0pt]
Faculty of Mathematics and Computer Science, \\[0pt]
Jagiellonian University \\[0pt]
ul. \L ojasiewicza 6, 30-348 Krak\'{o}w, Poland \\[0pt]
e-mail: Zbigniew.Hajto@uj.edu.pl\\[0pt]}

{\normalsize \vspace{0.5cm} ROUZBEH MOHSENI \\[0pt]
Faculty of Mathematics and Computer Science, \\[0pt]
Jagiellonian University \\[0pt]
ul. \L ojasiewicza 6, 30-348 Krak\'{o}w, Poland \\[0pt]
e-mail: Rouzbeh.Mohseni@doctoral.uj.edu.pl \\[0pt]}
 }

\end{center}

\vspace{0.2cm}

\begin{abstract}
In this paper we study the general concept of integrability in the broad
sense within the frame of differential Galois theory.  We concentrate on the
gradient systems which are not integrable. In spite of it, if we consider them as
the real dynamical systems,  they have trajectories with finiteness
properties of o-minimal type.
\end{abstract}

\vspace{0.2cm}

\section{Introduction}

From the time of publication of fundamental papers by Sergey Lvovich Ziglin
\cite{16,17}, the interest in the relation of Galoisian approach to
non-integrability and complex behaviour of dynamical systems has generated a
series of studies on the algebraic aspects of ``chaos theory". In the second
half of the XX century several interesting papers were published by R.
Churchill, J.J. Morales-Ruiz, J.M. Peris, D.L. Rod and B.D. Sleeman \cite{3a}%
, \cite{8}, \cite{13}. Quite recently K. Yagasaki and S. Yamanaka studied
Hamiltonian systems with two degrees of freedom and established a deep
connection between the existence of transverse heteroclinic orbits and the
non-integrability \cite{14,15}. In general, the relation of
Galoisian obstruction to integrability and chaotic dynamics is quite
complicated. In this paper we give a curious example of a non-integrable
gradient system which is not only non-chaotic but also is such that over $%
\mathbb{R}$ all its trajectories have tame topology in the sense of
o-minimal geometry (cf. \cite{4} and \cite{12}).

\section{Preliminaries}

\subsection{Integrability in the broad sense}

In this section we follow the notations of V.I. Arnold \cite{Ar}. Let $\Dot{x}=X_H$ be a Hamiltonian system and $F:M^{2n} \longrightarrow
\mathbb{R}$ be a first integral of it (i.e., $\{ H,F \} =0 $). It is possible
to choose canonical coordinates $x_1,...,x_n,y_1,...,y_n$ in a neighborhood
of a point $x_0 \in M$ with $d F(x_0)\neq 0$, such that $F(x,y)=y_1$. The
Hamiltonian function $H$ does not depend on $x_1$ in these coordinates, therefore putting $%
F=y_1=const$ gives us for $j\geq 2$

\begin{align}
\Dot{x_j}=\frac{\partial H}{\partial y_j},\,\,\,\,\, \Dot{y_j}=-\frac{%
\partial H}{\partial x_j}
\end{align}
which is a Hamiltonian system with $n-1$ degrees of freedom.

Therefore, having one first integral enabled us to lower the degree of
freedom of the system by one. Now to generalize this, assume that we have $%
n-1$ first integrals that are in involution, it can be shown that the degree
of freedom is reduced to one and the system will be ``integrable by
quadratures". This generalization can be stated as follows

\begin{theorem}
Let $M$ be a $2n$-dimensional symplectic manifold. A Hamiltonian vector
field $X_H$ on $M$ is integrable in the sense of Liouville, if there exist n
first integrals $F_j:M\longrightarrow \mathbb{R},\, 1\leq j\leq n$, such
that at every point of $M$ satisfy the following
\begin{enumerate}[(1)]
\item $dF_1,...,dF_n$ are linearly independent over a dense open set $U\subset M$;

\item  $\{ F_i,F_j \} =0$ for all $i,j$.
\end{enumerate}

\end{theorem}

The involutivity condition (i.e., the second property) implies that the
Hamiltonian system has an abelian Lie algebra $\mathcal{F}_a$ of
functionally independent first integrals $F_i$ with respect to the Poisson
brackets, and an abelian Lie algebra $\mathcal{S}_a$ of symmetries which
preserve the first integrals. O. Bogoyavlenskij \cite{2}
observed that these first integrals may be non-involutive and the symmetries
may be non-symplectic, therefore, the two Lie algebras are no longer
identified and gave the following definition.

\begin{definition}
The dynamical system $\Dot{x}=X$ on a smooth manifold $M^n$ is called
integrable in the broad sense if there exist
\begin{enumerate}[(1)]
\item $k$ functionally independent first integrals $F_1,...,F_k,\,\,1\leq k<n$,

\item an abelian ($n-k$)-dimensional Lie algebra $\mathcal{S}_a$ of symmetries
(i.e., vector fields $Y_i$ such that $Y_i(F_j)=0,\,\, \forall i,j.$), that are
linearly independent at the every point of $M$.

\end{enumerate}
\end{definition}

In Liouville integrability the dynamical system is Hamiltonian and these
two conditions are unified in the single condition of involutivity of the
first integrals. In general, a Hamiltonian system that is integrable in the
Liouville sense is also integrable in the broad sense but the converse is
not always true, therefore integrability in the broad sense is a
generalization of the Liouville integrability. A. Maciejewski and M.
Przybylska in \cite{7} called
the integrability in the broad sense, ``B-integrability".

A dynamical system on an n-dimensional manifold $N$ using the cotangent lift
can be included in the 2n-dimentional cotangent bundle $T^{*}N$. Let $%
\phi:N\longrightarrow N$ be a diffeomorphism, then we have the following

\begin{equation*}
\begin{tikzcd}
  T^*N \arrow[r,"\tilde{\phi}"] \arrow[d,"\pi "]
    & T^*N \arrow[d,"\pi "]  \\
    N  \arrow[r, "\phi "]
    &  N
 \end{tikzcd}
\end{equation*}
where $\tilde{\phi}$ is called the cotangent lift of $\phi$. If $X$ is a
vector field on the manifold $N$, using $X$ one can define the dynamical
system

\begin{align}
\Dot{x_i}=X_i,\, i=1,..., n.  \label{X}
\end{align}
Furthermore, let us suppose a function $f:T^*N\longrightarrow \mathbb{R}$
defined as follows

\begin{align}
f(p):=\langle p,X(x)\rangle,\, p\in T^*_xN.  \label{lift}
\end{align}
it is easy to show that the Hamiltonian vector field $X_f$ associated to $f$
is the following

\begin{align}
X_f=\left(X^i,-\frac{\partial f}{\partial x_i}\right).
\end{align}
and we say that $X_f$ is the cotangent lift of $X$. Let $X$ be a vector
field on $N$ and its cotangent lift be denoted by $Y$.

\begin{theorem}
If the vector field $X$ is (meromorphically) integrable in the broad sense,
then $Y$ is (meromorphically) Liouville integrable.

\begin{proof}
See Ayoul and Zung \cite{1}.
\end{proof}
\end{theorem}

In \cite{16,17} S.L. Ziglin started a new approach to integrability
using the monodromy group of the variational equation. Later Morales and
Ramis \cite{9,10}, took another approach to the integrability using
the differential Galois group of the variational equation and further
developed the theory. Morales and Ramis also stated that all the previous
Galoisian approaches to the Ziglin theory had a common
inconvenience, which is they restrict to the case where the variational equations belong to the
Fuchsian class (i.e., their singularities must be regular) and this
restriction is no longer necessary (for more details see the excellent survey \cite{3}). We now state the theorem by Morales and
Ramis which has attracted many researchers to the subject.

\begin{theorem}
\label{theorem 4} Assume that a complex analytical Hamiltonian system is
meromorphically completely integrable in a neighborhood of the integral
curve $x=x(t)$ (i.e., the solution of the Hamiltonian system). Then the
identity component of the differential Galois group of the variational
equations along the integral curve is Abelian.
\end{theorem}

As a remark, if the variational equation has irregular singularity at
infinity, then Theorem \ref{theorem 4} only provides obstruction to the
existence of rational first integrals.

Later in \cite{11} Morales, Ramis and Sim\'o proved this result for higher
variational equations. In \cite{1} Ayoul and Zung proved that the
differential Galois group of the variational equations (of any order) of the
original system is isomorphic to its counterpart for the cotangent lifted
system. Therefore, the condition that the system is integrable in Theorem %
\ref{theorem 4} can be generalized to being integrable in the broad sense
and the group being the differential Galois group of variational equations
of any order (see also \cite{7} Theorem 13).

\begin{theorem}
\label{theorem 5} Assume that a dynamical system is meromorphically
integrable in the broad sense in a neighborhood of the integral curve $x=x(t)
$ (i.e., the solution of the dynamical system). Then the identity component
of the differential Galois group of the variational equations (of any order)
along the integral curve is Abelian.
\end{theorem}

\subsection{Pfaffian manifolds}

In this section we gather together some facts from the Khovanskii's theory
of fewnomials \cite{5}, which we shall use to the study of gradient
trajectories. We assume that the ambient space is $\mathbb{R}^n$ and we
consider real analytic foliations of codimension 1 (cf. \cite{CH} Section
I.4). For the convenience of the reader, let us recall the definition of a
Pfaffian hypersurface and its topological properties.

\begin{definition}
A Pfaffian hypersurface of $\mathbb{R}^n$ is a triple $(V,\mathcal{F},M)$,
where

\begin{enumerate}[(1)]
\item $M$ is an open semianalytic set of $\R^n$,
\item $\F$ is a foliation of codimension 1 defined in an open neighbourhood of the closure of $M$,
\item $V$ is a leaf of $\F_M$, i.e., is a connected maximal integral submanifold of the restricted foliation,
\item $\Sing \F \cap M=\emptyset$.

\end{enumerate}
\end{definition}

A Pfaffian hypersurface is called \emph{separating} if $M\setminus V$ has
two connected components and $V$ is their common boundary in $M$. A Pfaffian
hypersurface has the \emph{Rolle property} (i.e., is \emph{Rollian}) if each
analytic path $\gamma:[0,1] \rightarrow M$, such that $\gamma(0), \gamma(1)$
are points in $V$, has at least one point $\gamma(t)$ such that the tangent
vector $\gamma^{\prime}(t)$ is tangent at the foliation $\mathcal{F}$ (i.e.,
if at this point $\mathcal{F}$ is determined by a local generator $\omega$
then $\gamma^{\prime}(t) \in \Ker \omega$). By the theorem of
Khovanskii-Rolle, a separating hypersurface is Rollian but the converse is
not true.

In the context of Pfaffian geometry, Khovanskii's finiteness theorem has the
following form (cf. \cite{CH} Theorem 1).

\begin{theorem}
Let $M$ be an open semianalytic set in $\mathbb{R}^n$ and $X\subset M$
semianalytic and bounded in $\mathbb{R}^n$. For each finite collection of
Pfaffian hypersurfaces $(V_1,\mathcal{F}_1,M),\dots,$ \linebreak $(V_p,%
\mathcal{F}_p,M)$ which have the Rolle property for the paths in $X$, there
exists a number $b_0 \in \mathbb{N}$ (depending only on $M, X, \mathcal{F}%
_1,\dots,\mathcal{F}_p$) such that the number of connected components of $%
X\cap V_1\cap \dots \cap V_p$ is smaller than $b_0$.
\end{theorem}

\begin{remark}

Separating Pfaffian hypersurfaces are interesting examples of ``tame spaces" as
predicted by Grothendieck in his sketch of a programme \cite{12}. Separating Pfaffian hypersurfaces in the real plane $\mathbb{R}^2$ are called the \emph{Pfaffian curves}.
Pfaffian curves have similar properties to the semianalytic arcs, therefore, their
topology can be described within the theory of o-minimal structures (cf.
\cite{4},\cite{Sp}).
\end{remark}

\section{Gradient trajectories in the plane}

\subsection{Non-integrability of planar vector fields}

In a very interesting paper \cite{Ac}, P.B. Acosta-Hum\'{a}nez, J.T. L\'{a}%
zaro, J.J. Morales-Ruiz and C.Pantazi tackled non-integrability of
complex planar vector fields of the form $X=P\frac{\partial }{\partial x}+Q%
\frac{\partial }{\partial y}$ with $P,Q\in
\mathbb{C} \lbrack x,y]$ (cf. \cite{Ac} Section I, Theorem B). In fact they have
elaborated an effective algorithm in order to check the necessary conditions for the integrability
of polynomial vector fields in the plane. Their method is illustrated with several families of examples,
however, in the presented families there are no examples with scalar potentials, i.e., gradient vector fields.
On the other hand, gradient trajectories are intensively studied in the range of problems related with the
stability of critical points of the gradient systems (see e.g. \cite{6} Chapter I, Section 8).
Let us note the common conviction that gradient trajectories have quite ``moderated" topology.
More precisely, K.Kurdyka, T.Mostowski and A. Parusi\'nski in their landmark paper \cite{KK}, proved the
gradient conjecture of R. Thom. Simultaneously, they have stated an even more general conjecture called
the finiteness conjecture for the gradient systems. This conjecture in the case of plane
is easily proved by the elementary methods of subanalytic geometry (see \cite{KK}, Proposition 2.1).
Let us consider the following gradient dynamical system

\begin{equation}\label{eq5}
\dot{x}= \grad F, \text{\ \ where \ } F \in \R[x,y].
\end{equation}

\noindent
Taking into account that for $n$ = 2 (i.e., in the real plane) any subanalytic set is actually semianalytic, it is easy to see that gradient trajectories of system (\ref{eq5})  are finite unions of Pfaffian curves (cf. Remark 2.2.3).

\begin{example}\label{example}
Consider the planar gradient field

$$
\dot{x}= \grad F, \text{\ \ where \ } F(x,y)= \dfrac 1 3 x^3+\dfrac 1 2 x^2+(x+y)^2y^2+\dfrac 1 4 y^4
$$

\noindent
and its associated foliation

\begin{equation}\label{eq6}
\dfrac{dy}{dx} = \dfrac {\partial F / \partial y}{\partial F / \partial x}=\dfrac{2(x+y)y^2+2(x+y)^2y+y^3}{x^2+x+2(x+y)y^2}
\end{equation}

\noindent
Following \cite{Ac} Section 3, we will show that the linearized second variational equation $LVE_2$ of (\ref{eq6}) along a particular solution has non-commutative identity component $G_2^0$ of its Galois group $G_2$.

First let us observe that the straight line $\Gamma =\{ y=0 \}$ is an invariant curve of (\ref{eq6}) and we have the second linearized variational equation along $\Gamma$

\begin{equation}\label{eq7}
\begin{array}{lcl} \chi_1'&=& 2 \dfrac {2x}{x+1} \chi_1 \\[16pt]
 \chi_2'&=&  \dfrac {2x}{x+1} \chi_2+ \dfrac{12}{x+1} \chi_1
\end{array}
\end{equation}

\noindent
In the notations of \cite{Ac} Section 3 we have

$$
\begin{array}{c}
\beta_1(x)= \dfrac{2x}{x+1} , \, \beta_2(x)= \dfrac{12}{x+1}, \\[16pt]
\omega= \dfrac{e^{2x}}{(x+1)^2}, \, \theta_1= \displaystyle \int \dfrac{12 e^{2x}}{(x+1)^3}
\end{array}
$$

\noindent
In the notation of Maple 2019, we obtain

\begin{equation}\label{eq8}
\theta_1= -\dfrac{(12x+18)e^{2x}}{(x+1)^2}-24e^{-2}Ei_1(-2x-2),
\end{equation}

\noindent
where $Ei_1(z)$ denotes the special function \textquotedblleft exponential integral \textquotedblright.
In order to apply Proposition 3.1 from \cite{Ac} Section 3, we have to check the hypotheses $(H_1)$ and $(H_2)$. $(H_1)$ is clearly fulfilled. Concerning $(H_2)$ let us suppose that $\theta_1$ is a rational function in $x$ and $\omega$. Then, since from the expression for $\theta_1$ in (\ref{eq8}), we obtain

$$Ei_1(-2x-2)=-\dfrac 1 {24e^{-2}} \left( \theta_1+ (12x+18)\omega \right),$$

\noindent
$Ei_1(-2x-2)$ would be rational in $x$ and $\omega$, a contradiction, since $Ei_1(-2x-2)$ is not an elementary function.
\end{example}

\section{Final remarks and conclusions}

As R.C. Churchill and D.L. Rod observed in their joint note \cite{3a}, an important consequence of chaos in a Hamiltonian system is its non-integrability. On the other hand the study of chaotic dynamical systems from the algebraic point of view is still in the infancy and as shows e.g Example \ref{example}, non-integrability in the Arnold-Liouvillian sense may be unrelated with the complexity of the topology of solutions (cf \cite{5}, Introduction). The Morales-Ramis theory in its classical form \cite{9}, \cite{10} deals with the problems of non-integrability in the complex case mainly because it is based on the Picard-Vessiot theory with the additional assumption that the field of constants is algebraically closed. However from the time of discovery of so called \textquotedblleft real Picard-Vessiot fields \textquotedblright \cite{CHP} it seems very interesting to develop an analogous theory over the field of real numbers. In such a case it would be possible to study \textquotedblleft tame topology \textquotedblright of the solutions following the theory of fewnomials. In fact in the case of linear systems the relation of real Liouville extensions with Khovanskii-Gel'fond theory is presented in \cite{CH}.

\section*{Acknowledgments} The authors especially thank the participants of the Differential Galois Theory Seminar at Jagiellonian University where many topics of our paper were discussed. The first author thanks Teresa Crespo and Askold Khovanskii for valuable comments on real Liouville integrability. The second author would like to thank Alexey Bolsinov for valuable comments concerning the theory of chaos.

\end{document}